\documentclass[11pt]{article}
\usepackage{amsfonts}
\usepackage{amssymb}
\usepackage{amsfonts,epsfig}
\begin{document}

\newtheorem{them}{Theorem}
\newtheorem{lemm}{Lemma}
\newtheorem{prop}{Proposition}
\newtheorem{definition}{Definition}
\newtheorem{rema}{Remark}
\newtheorem{coro}{Corollary}
\newtheorem{exam}{Example}

\newcommand{\ssq}{$\square$}
\newcommand{\nus}{\smallsetminus}
\newcommand{\cl}{{\cal L}}
\newcommand{\qed}{\hfill$\Box$}
\newcommand{\pr}{{\bf Proof. }}
\newcommand{\rr}{\Bbb{R}}
\newcommand{\cc}{\Bbb{C}}
\newcommand{\wh}{\widehat}
\newcommand{\pl}{\partial}
\newcommand{\ol}{\overline}
\newcommand{\ul}{\underline}
\newcommand{\ra}{\Rightarrow}
\newcommand{\beq}{\begin{equation}}
\newcommand{\eeq}{\end{equation}}
\newcommand{\beqn}{\begin{eqnarray}}
\newcommand{\eeqn}{\end{eqnarray}}
\newcommand{\beqnn}{\begin{eqnarray*}}
\newcommand{\eeqnn}{\end{eqnarray*}}
\newcommand{\vep}{\varepsilon}
\newcommand{\vno}{\varnothing}
\newcommand{\al}{\alpha}
\newcommand{\la}{\lambda}
\newcommand{\Lt}{{\cal L}(A,\alpha)}
\newcommand{\rM}{{\cal M}_{n}(\mathbb{R})}
\newcommand{\cM}{{\cal M}_{n}(\mathbb{C})}
\newcommand{\M}{{\cal M}_{n}}

\title{ Further results on the Craig-Sakamoto Equation }

\author { John Maroulas$^{1}$}
\maketitle

\vspace{0.2cm}

\footnotetext[1]{Department of Mathematics, National Technical
University, Zografou Campus, Athens 15780, GREECE.
E-mail:maroulas@math.ntua.gr.  This work is supported by a grant
of the EPEAEK, project "Pythagoras II". }

\vspace{0.2cm}

\vspace{0.2cm}

\begin{abstract}
In this paper necessary and sufficient conditions are stated  for
the Craig-Sakamoto equation
 $\,det(I-sA-tB)=det(I-sA)\,det(I-tB),\,$  for all
scalars $\,s,\,t.\,$ Moreover, spectral  properties for $\,A\,$
and $\,B\,$ are investigated.
\end{abstract}
\ \\

\vspace{0.15cm}


\noindent
\section{Introduction}

~~~~~Let $\,M_n(\mathbb{C})\,$ be the set of $\,n \times n\,$
matrices with elements in $\,\mathbb{C}.\,$  For $\, A\,$ and $\,B
\in M_{n}(\mathbb{C}),\,$ the well known in Statisticts \cite{Du}
Craig-Sakamoto (CS) equation
\begin{eqnarray}
\label{eq1} det(I-sA-tB)=det(I-sA)\,det(I-tB) \end{eqnarray} for
all scalars $\,s,\,t\,$  has occupied several researchers. In
particular, in \cite{T} O. Trussky  presented that the CS equation
is equivalent to $\,AB=O,\,$ when $\,A,\,B\,$ are normal and most
recently in \cite{O} Olkin  and in \cite{L} Li  proved the same
result in a different way. The author, together with M. Tsatsomero
and P. Psarrako in \cite{MPT}, have investigated the CS equation
involving the eigenspaces  of $\,A,B\,$ and $\,sA + tB .\,$ Being
more specific, if $\,\sigma(X)\,$ denotes the spectrum for a
matrix $\,X,\,$ $\, m_X (\lambda) \,$
 the algebraic multiplicity of $\, \lambda \in \sigma(X),\,$ and
 $\,E_X(\lambda)\,$ the generalized eigenspace corresponding to
 $\,\lambda,\,$ we have shown in \cite{MPT}:
\newpage
\begin{prop}
\label{prop1} ~ For the $\, n \times n\,$ matrices  $\,A,\,B\,$
the following are equivalent :
\begin{description}
    \item[ I. ]   ~~~The CS equation holds
    \item[ II. ] ~~for every $\,s,\,t \in \mathbb{C,}\,$ $\, \sigma ( sA
    \oplus tB ) = \sigma \left( ( sa + tB )  \oplus O_n
    \right), \,$ where $\,O_n\,$ denotes the zero matrix
    \item[ III. ] ~$ \sigma ( sA
    + tB ) = \left\{ \, s \mu_i + t \nu_i \,:\, \mu_i \in \sigma (A) ,\;\,\, \nu_i \in \sigma(B )\,\right\},
    \,$ where the pairing of eigenvalues ~~~~~~~~~~~~~requires  either $\,\mu_i
    = 0\,$ or $\,\nu_i
    = 0.\,$
\end{description}
\end{prop}
\begin{prop}
\label{prop2} ~ Let the $\, n \times n\,$ matrices $\,A,\,B\,$
satisfy the CS equation. Then,
\begin{description}
    \item[ I. ]  ~~~$\,m_A (0) + m_B (0)\, \geq \, n.\,$
    \item[ II. ]  ~~If $\,A \,$ is nonsingular, then $\,B\,$ must be nilpotent.
    \item[ III. ]   ~If $\,\lambda = 0\,$ is semisimple eigenvalue of $\,
    A\,$ and $\,B,\,$ then  $\,rank(A) + rank(B) \leq n .\,$
\end{description}
 \end{prop}
\begin{prop}
\label{prop3} ~ Let $\, \lambda =0\,$ be semisimple eigenvalue of
$\, n \times n\,$ matrices $\,A\,$ and $\,B\,$ such that
$\,BE_{A}(0) \subset E_{A}(0).\,$ Then the following are
equivalent.
\begin{description}
    \item[ I. ]  ~~~~Condition CS holds.
    \item[ II. ]  ~~$\,\mathbb{C}^{n}=E_{A}(0)+E_{B}(0). \,$
    \item[ III. ]   ~$\,AB=O.\,$
\end{description}
 \end{prop}
  The remaining results in \cite{MPT} are based on the
 basic assumption that $\,\lambda = 0\,$ is a semisimple
 eigenvalue of $\,A\,$ and $\,B.\,$ Relaxing this restriction, we
 shall attempt here to look at the CS equation focused on the factorization of polynomial of two variables $\,f(s,\,t)= det(I-sA-tB).$
\newline Also, considering the determinants in (\ref{eq1}), new conditions necessary and sufficient on CS property are stated.

\noindent
\section{Spectral results }

~~~~~The first statement on the CS property is obtained
investigating the determinantal equation through the Theory of
Polynomials. By Proposition \ref{prop2} {\bf II}, it is clear that
the CS equation is worth valuable when the $\,n \times n\,$
matrices $\,A\,$ and $\,B\,$ are singular. Especially, we define
that \newline "{\it $\,A\,$ and $\,B\,$ are called $\,r$-{\bf
complementary}, if and only if  {\rm at most}, $\,r\,$ rows
(columns), $\,a_{i_{1}}, \,a_{i_{2}},\, \cdots, a_{i_{r}}\,$ of
$\,A\,$ are shifted and substituted by the corresponding
$\,b_{i_{1}}, \,b_{i_{2}},\, \cdots, b_{i_{r}}\,$ rows (columns)
of $\,B,\,$ such that the structured matrix $\,N(i_{1}, \,i_{2},\,
\cdots, i_{r})\,$ of $\,a$' s and $\,b$' s rows is nonsingular}."
\vspace{0.3cm} \newline Note that, $\, n-r \leq rank(B).\,$

For example, the pair of matrices
$$\,A=\left[
\begin{array}{ccc}
  0 & 1 & 0 \\
  0 & 0 & 1 \\
  0 & 0 & 0 \\
\end{array}%
\right] ,\;\;\;\;\,B=\left[
\begin{array}{ccc}
  1 & 0 & 0 \\
  0 & 0 & 1 \\
  0 & 0 & 0 \\
\end{array}%
\right] \,$$ is not $\,1\,$ or $\,2-$complementary, on behalf of
$\,rank \left[
\begin{array}{c}
  A \\
  B \\
\end{array}%
\right] =3, \,$ but the pair
$$\,{\cal A}=\left[
\begin{array}{ccc}
  0 & 0 & 0 \\
  0 & 1 & 1 \\
  0 & 0 & 1 \\
\end{array}%
\right] ,\;\;\;\;\,{\cal B} = B\,$$ is $\,1$-complementary and not
$\,2$-complementary, since $\,detN(b_1,a_2,a_3)=det\left[
\begin{array}{ccc}
  1 & 0 & 0 \\
  0 & 1 & 1 \\
  0 & 0 & 1 \\
\end{array}%
\right] \neq 0\,\,$ and $\,detN(b_1,b_2,a_3)=det\left[
\begin{array}{ccc}
  1 & 0 & 0 \\
  0 & 0 & 1 \\
  0 & 0 & 1 \\
\end{array}%
\right] = 0.$

\vspace{0.5cm}

\begin{prop}
\label{p1} ~Let the $\,n \times n\,$ singular matrices $\,A\,$ and
$\,B\,$ be $\,[n-m_B(0)]$-complementary with
$\,\theta=\displaystyle{\sum_{i_{1},\ldots,i_{n-m_{B}(0)}}detN(i_{1},i_{2},\ldots,i_{n-m_{B}(0)})}
\neq 0,\,$ where the sum is over all possible combinations
$\,i_{1},\ldots,i_{n-m_{B}(0)}\,$ of $\,n-m_{B}(0)\,$ of the
indices $\, 1,2,\ldots, n.\,$ If they satisfy the CS equation,
then
$$\,m_A (0) + m_B (0)\, = \, n.\,$$
\end{prop}
\pr~ ~Let $\, rank B = b\, ( < n). \,$ Then $\,\lambda =0 \,$ is
eigenvalue of $\, B\,$ with algebraic multiplicity $\, m_B(0)=m
\geq n-b. \,$ Denoting
$$
\beta(t) \doteqdot  det( t I -B) = t^{n} + \beta_1 t^{n-1} +
\cdots + \beta_{n-m} t^{m},
$$
where $\, \beta_k=(-1)^{k}\sum B_{k} \,$ and $\, B_{k}\,$ are the
$\, k \times k\,$ principal minors of $\,B,\,$ then
\begin{eqnarray*}
det( t B -I ) & = & (-1)^n \, t^{n} \, det ( t^{-1} I -B )\\
& = & (-1)^n \, \left( 1 + \beta_{1}\, t + \cdots + \beta_{n-m}
\,t^{n-m}  \right ).
\end{eqnarray*}
The polynomial $\,\widetilde{\beta}(t)=1+\beta_{1}\, t + \cdots
+\beta_{n-m} \,t^{n-m}\,$ has precisely $\, n-m\,$ nonzero roots,
let $\,t_1, t_2, \cdots, t_{n-m},\,$ since $\,\widetilde{\beta}(0)
=1 \neq 0.\,$ Moreover, we have
\begin{eqnarray}
\label{eq2} det(sA + tB-I)  =  |A| s^{n} +f_1(t)s^{n-1} + \cdots +
f_{n-1}(t) s + |tB-I| ,
\end{eqnarray}
where
\begin{eqnarray*}
f_1(t) = \sum _{i} det \hat{A}_{i},\;\; \mbox {  with  } \;\;
\hat{A}_{i}=\left[
 \begin{array}{ccccc}
 a_{11} & & \cdots & & a_{1n} \\
\vdots & &  & & \vdots \\
tb_{i1} & \cdots & t b_{ii}-1 &  \cdots  & t b_{in} \\
\vdots & &  &  & \vdots \\
a_{n1} & & \cdots & & a_{nn} \\
\end{array}
 \right].
 \end{eqnarray*}
Note that, $\,\widehat{A}_i\,$ arises by $\,A\,$ when the
$\,i-$row of $\,A\,$ is substituted by the $\,i-$row of
$\,tB-I.\,$ Similarly,
\begin{eqnarray*}
f_2(t) = \sum _{i,j} det \hat{A}_{ij},\;\; \mbox {  with  }
\;\;\hat{A}_{ij}=\left[
 \begin{array}{cccccc}
 a_{11} &  & \cdots & & & a_{1n} \\
 \vdots & & & & & \vdots \\
t b_{i1} & \cdots & t b_{ii}-1 & \cdots & &  t b_{in} \\
\vdots & & \ddots &  & & \vdots \\
t b_{j1} & \cdots & &  t b_{jj}-1 &  \cdots & t b_{jn} \\
 \vdots & & & & & \vdots \\
a_{n1} & &  \cdots & & &  a_{nn} \\
\end{array}
 \right],
  \end{eqnarray*}
and $\,\hat{A}_{ij}\,$ is obtained by $\,A,\,$ substituting the
$\,i\,$ and $\,j\,$ rows of $\,A\,$ by the corresponding rows of
$\,tB-I.\,$ The summation in $\,f_2(t)\,$ is referred to all pairs
of indices $\,i,j\,$ by $\,\{1,2,\ldots, n\}.\,$ Hence, by the
equation (\ref{eq2}) and the CS equation
$$
(-1)^{n} \, det(sA+tB-I)=det(sA-I)\, det(tB-I), \;\;\;\;\;\;\;
\forall  \; s,\,t
$$
for $\, t=t_1,\,t_2,\, \cdots , t_{n-m},\,$ we obtain
\begin{eqnarray*}
|A| s^{n} +f_1(t_{i})s^{n-1} + \cdots + f_{n-1}(t_{i}) s =0,
\;\;\;\;\;\; \forall \;\; s
\end{eqnarray*}
and consequently
\begin{eqnarray}
\label{eq3} |A|=0\,\;\;, \;\;\, f_1(t_{i})=f_2(t_{i})= \cdots
=f_{n-1}(t_{i})=0,\,\;\;\mbox{ for } \;\;\,\,i=1,2,\ldots, n-m\,.
\end{eqnarray}
Due to the matrices $\,A\,$ and $\,B\,$ are
$\,[n-m_B(0)]$-complementary and the leading coefficient of
$\,f_{n-m}(t)\,$ is equal to the nonzero $\,\theta,\,$ then
$\,\,deg( f_{n-m}(t) )= n-m \,$ and $\,\,deg( f_{k}(t) ) \leq n-m
,\,$ for  $\,\,k=1,2,\ldots, n-m-1.\,$ Moreover, by (\ref{eq3}) we
have
$$
f_1(t)= f_2(t)= \cdots = f_{n-m-1}(t)= 0, \;\;\; \forall \;\; t
$$
Reminding that $\,A_{\ell}\,$ denotes the $\, \ell \times \ell\,$
principal minor of $\,A,\,$ by $\,f_1(t) = 0,\,$ clearly
\begin{eqnarray*}
f_1(0)= \sum A_{n - 1} = 0  \;\; \Longrightarrow  \;\; c_{n-1} =
0.
\end{eqnarray*}
Similarly, by
\begin{eqnarray*}
& & f_2(t) =  0\,\Longrightarrow \;\; \sum A_{n-2} = 0 \;\;\;
\Longrightarrow \;\;\;c_{n-2}=0
\\
& & \;\; \;\; \vdots
\\
& & f_{ n - m -1 }(t) = 0 \;\; \Longrightarrow \;\; \sum A_{m+1} =
0 \;\;\; \Longrightarrow \;\;\;c_{m+1}=0 ,
\end{eqnarray*}
and consequently
\begin{eqnarray}
\label{eq4} \delta_A(\lambda) = |\lambda I - A| & = &  \lambda
^{n} - c_1 \lambda
^{n-1} + c_2 \lambda ^{n-2} + \cdots + (-1) ^{n} |A|  \nonumber \\
& = &  \lambda ^{n} - c_1 \lambda ^{n-1} + c_2 \lambda ^{n-2} +
\cdots + (-1) ^{m } c_{ m } \lambda ^{ n - m }.
\end{eqnarray}
In (\ref{eq4}), $\,c_{m} \neq 0, \,$ since $\,(-1)^{n-m}c_m
=\theta t_1t_2\,\cdots\,t_{n-m}.\,$  Thus, $\, \lambda = 0 \,$ is
eigenvalue of $\,A\,$ with algebraic multiplicity $\, n -
m_{B}(0),\,$ whereby we conclude
$$\,m_A (0) + m_B (0)\, = \, n.\,$$ \qed

\begin{rema}
~ \rm{  By the proof of Proposition \ref{p1}, it is evident that
the equality  $\,\,m_A (0) + m_B (0)\, = \, n\,$ holds, when the
matrices  $\,B\,$ and $\,A\,$  are $\,[n-m_A(0)]$-complementary
and
$$
\,\theta=\sum_{j_{1},\ldots,j_{n-m_{A}(0)}}detN(j_{1},j_{2},\ldots,j_{n-m_{A}(0)})
\neq 0.\,$$ }
\end{rema}

\vspace{0,2cm}

\begin{coro}
 \label{c1}
~ Let the $\, n \times n \,$ singular and
$\,[n-m_B(0)]$-complementary matrices $\,A\,$ and $\,B.\,$ If
$\,\theta \neq 0\,$ and these matrices satisfy the CS equation
(\ref{eq1}), then
\begin{description}
    \item[ I. ]  ~~$\,\lambda = 0\,$ is semisimple eigenvalue of
    $\,A\,$ and $\,B\,\Longrightarrow \,rank A + rank B =  n .\,$
    \item[ II. ] ~$\,\lambda = 0\,$  is semisimple eigenvalue   of
    $\,A\,$  $\,\Longrightarrow rank A = m_B(0). \,\,$
\end{description}
\end{coro}
\pr ~{\bf I. }~ Because $$\, n-rank A \leq m_A(0) = n -
m_B(0),\,$$ we have $\, rank A + rank B \geq m_B(0) + r \geq n.\,$
Hence, by {\bf III}, Proposition \ref{prop2}, we obtain the
equality.
\newline ~{\bf II.}~
By the assumption and Proposition \ref{p1} we have $\, rank A = n-
m_A(0) = m_B(0).\, $ \qed

 \vspace{0,3cm}
Closing this section, we present a property of generalized
eigenspaces of nonzero eigenvalues of $\,A\,$ and $\,B.\,$

 \vspace{0,3cm}

\begin{prop}
\label{p2} ~ Let $\,\lambda = 0 \,$ be semisimple eigenvalue of
$\,n \times n\,$ matrices $\,A\,$ and $\,B\,$  such that $\,E_A(0)
+ E_B(0) = \mathbb{C}^{n}.\,$ If for any $\,\lambda \in \sigma(A)
\backslash \{0\}\,$, (or, $\,\mu \in \sigma(B) \backslash
\{0\}\,$), the corresponding generalized eigenspaces
$\,E_A(\lambda),\,$ ($\, E_B(\mu)\,$) satisfy $\,\,E_A(\lambda)
\subseteq  E_B(0),\,$ (or, $\,E_B(\mu) \subseteq E_A(0)\,$), then
\begin{description}
    \item[ I. ]   ~~~$\,A,\,B\,$ have the CS property.
    \item[ II. ] ~~$\,E_{A}(\lambda)= E_{I-sA-tB}(1-s\lambda),\,$ and $\,E_{B}(\mu)= E_{I-sA-tB}(1-t\mu).\,$
\end{description}
\end{prop}
\pr~{\bf I.}~Since $\,E_A(\lambda) \subseteq E_B(0),\,$ for every
$\,w=w_1 + w_2 \in \mathbb{C}^{n},\,$ where $\,w_1 \in
\bigoplus_{\lambda}E_A(\lambda),\,$ $\,w_2 \in E_A(0),\,$ we have
$\,BAw=BA(w_1+w_2)=BAw_1=0.\,$ Thus, $\,BA=O\,$ and consequently
$\,AE_B(0) \subseteq E_A(0).\,$ The assumption $\,E_A(0)+
E_B(0)=\mathbb{C}^{n},\,$ and Proposition \ref{prop3}, lead to the
statement {\bf I}.
\newline{\bf II.}~ Let $\,\lambda \in \sigma(A)
\backslash \{0\},\,$ and $\,x_k\in E_A(\lambda)\,$ be generalized
eigenvector of $\,A\,$ of order $\,k.\,$ By the assumption,
$\,x_k\in E_B(0),\,$ and yields
\begin{eqnarray*}
(I - sA -tB)x_k & = & (I - sA )x_k = x_k -s (\,\lambda x_k +
x_{k-1}\,) \\
& = &  ( 1-s \lambda ) x_k -s x_{k-1}.\,
\end{eqnarray*}
Thus, for all chain $\,x_1,\,\ldots,\, x_k, \ldots,\,x_{\tau}\,$
of $\,\lambda,\,$ we have
\begin{eqnarray}
\label{5} (I - s A - t B )
\left[%
\begin{array}{ccc}
  x_1 & \ldots & x_{\tau} \\
\end{array}%
\right] =
\left[%
\begin{array}{ccc}
  x_1 & \ldots & x_{\tau} \\
\end{array}%
\right]
\left[%
\begin{array}{cccccc}
  1- s\lambda  & -s  & & &  \\
    0 & 1- s\lambda   & \;\;  -s & & O &   \\
 \vdots &  & \ddots  &   & \ddots  &  \\
  & & &  &  1- s \lambda  & -s  \\
      0 & 0   &  & & &  1- s \lambda
\end{array}%
\right]_{\tau \times \tau }
\end{eqnarray}
Moreover, by the statement {\bf III} in Proposition \ref{prop1},
$\,s\lambda\,$ and $\,t\mu \in \sigma(sA+tB).\,$ The equivalence
of CS equation and $\,\mathbb{C}^{n}=E_A(0)+ E_B(0)\,$ in
Proposition \ref{prop3} and the assumption $\,E_A(\lambda)
\subseteq E_B(0), \,$ lead to $\,E_B(\mu) \subseteq E_A(0).\,$
Similarly, if $\,y_{\ell} \in E_B(\mu)\,$ is generalized
eigenvector of order $\,\ell,\,$ then $\,y_{\ell} \in E_A(0)\,$
and
\begin{eqnarray*}
 (I - sA -tB)y_{\ell} & = & (I - tB )y_{\ell} = y_{\ell} - t
(\,\mu y_{\ell} +
y_{\ell -1}\,) \\
& = &  ( 1-t \mu ) y_{\ell} -t y_{\ell-1},\,
\end{eqnarray*}
and for all chain $\,y_1,\ldots, \,y_{\ell},
\ldots,\,y_{\sigma}\,$ we obtain
\begin{eqnarray}
\label{6} (I - s A - t B )
\left[%
\begin{array}{ccc}
  y_1 & \ldots & y_{\sigma} \\
\end{array}%
\right] =
\left[%
\begin{array}{ccc}
  y_1 & \ldots & y_{\sigma} \\
\end{array}%
\right]
\left[%
\begin{array}{cccccc}
  1- t\mu & -t  & & &  \\
    0 & 1- t \mu  & \;\;  -t & & O &   \\
 \vdots &  & \ddots  &   & \ddots  &  \\
  & & &  &  1- t \mu & -t  \\
      0 & 0   &  & & &  1- t \mu
\end{array}%
\right]_{\sigma \times \sigma }
\end{eqnarray}
Clearly, by (\ref{5}) and (\ref{6}) are implied the equations in
{\bf II}, for any $\,s,\,t.\,$ \qed

\begin{rema}
~ \rm{ For $\, z \in E_A(0) \bigcap E_B(0)\,$ obviously $ \,( I -
sA - t B ) z= z ,\,$ $ \, \forall \, s,\, t.\,$ Therefore by the
above proposition the Jordan canonical form of $\,I - sA - t B,\,$
and the matrix
$$
F= I_{\nu} \bigoplus _{\lambda_{A} \neq 0}
\left[%
\begin{array}{cccc}
  1- s \lambda_{A}  & -s  & & O \\
    & 1- s\lambda_{A}   & \ddots &     \\
  &   & \ddots   & -s   \\
  O &  &   & 1- s \lambda_{A}
\end{array}%
\right] \bigoplus _{\mu_{B} \neq 0}
\left[%
\begin{array}{cccc}
  1- t \mu_{B}  & -t  & & O \\
    & 1- t \mu_{B}  & \ddots &     \\
  &   & \ddots   & -t   \\
  O &  &   & 1-t \mu_{B}
\end{array}%
\right],
$$
are similar.

The order $\,\nu \,$ of submatrix $\,I_{\nu}\,$ of $\,F\,$
 declares the number of linear independent eigenvectors which
correspond to the eigenvalue  $\,\lambda =1\, $ of $\,I-sA-tB.\,$
Clearly, theses eigenvectors belong to $\,\, E_B(0 ) \backslash
E_A(\lambda ), \,\;E_A(0 ) \backslash E_B(\mu ),\,$ and $\, E_A(0
) \bigcap E_B(0 ),\, $ and $\,\nu\,$ is equal to
$$
\nu = n - (rankA + rankB ) = n - \left( dim \bigcup_{\lambda \neq
0} E_A(\lambda) + dim \bigcup_{\mu \neq 0} E_B(\mu ) \right).
$$
}
\end{rema}
\vspace{0.3cm}
 \noindent
\section{Criteria for CS equation}

Let
\begin{eqnarray}
\label{eqc7} f(s,t)=det(I-sA-tB)= \sum_{p,q=0}^{n}
m_{pq}s^{p}t^{q},\;\;\;\;\;\;p+q \leq n.
\end{eqnarray}
Denoting by $\,
\;x=\left[%
\begin{array}{ccccc}
  1 & s & s^{2} & \cdots & s^{n} \\
\end{array}%
\right]^{T} , \;\;\;\;\;
y=\left[%
\begin{array}{ccccc}
  1 & t & t^{2} & \cdots & t^{n} \\
\end{array}%
\right]^{T} ,\,\;\, $ then (\ref{eqc7}) is written obviously
$$\,f(s,t) =x^{T}My,\,$$
where $\,M=\left[m_{pq} \right]^{n}_{p,q=0}\,$, with $m_{00}=1.\,$

\begin{prop}
\label{p3}~ Let $\,A,\,B \in M_n(\mathbb{C}).\,$ The CS equation
holds for the pair of matrices $\,A\,$ and $\,B\,$ if and only if
$\,\, rank M = 1.\,$
\end{prop}
\pr Let $A$ and $B$ managed by the CS property. Then the equation
(\ref{eq1}) is formulated as
\begin{eqnarray}
\label{eqc8} x^{T}My=x^{T}a \, b^{T}y,
\end{eqnarray}
where
$$
a=\left[%
\begin{array}{cccc}
  1 & a_{n-1} & \cdots & a_{0} \\
\end{array}%
\right]^{T} , \;\;\;\;\;\;
b=\left[%
\begin{array}{cccc}
  1 & b_{n-1} & \cdots & b_{0} \\
\end{array}%
\right]^{T} ,
$$
and $\,a_{i},\,b_{i}\,$ are the coefficients of characterictic
polynomials
$$
det(\lambda I -A)=\lambda^{n} + a_{n-1} \lambda^{n-1}+ \ldots +
a_{0},\;\;\;det(\lambda I -B)=\lambda^{n} + b_{n-1} \lambda^{n-1}+
\ldots + b_{0}.
$$
Hence, by (\ref{eqc8}) for any $\,s_1 \neq s_2 \neq \cdots \neq
s_{n+1} \,$ and $\,t_1 \neq t_2 \neq \cdots \neq t_{n+1} \,$ we
have
\begin{eqnarray}
\label{eqc9} V^{T}\left( M- a \, b^{T}\right)W=O,
\end{eqnarray}
where
$$
V=\left[%
\begin{array}{ccc}
  1 & \cdots  & 1 \\
  s_1 & \cdots & s_{n+1} \\
  \vdots &  & \vdots \\
  s^{n}_1 & \cdots & s^{n}_{n+1} \\
\end{array}%
\right],\;\;\;\;\;\;W=\left[%
\begin{array}{ccc}
  1 & \cdots  & 1 \\
  t_1 & \cdots & t_{n+1} \\
  \vdots &  & \vdots \\
  t^{n}_1 & \cdots & t^{n}_{n+1} \\
\end{array}%
\right].
$$
Clearly, by (\ref{eqc9}), we recognize that $\,M= a \,b^{T},\,$
i.e., $\,rankM=1.$

Conversely, if $\,rankM=1, \,$ then $\,M=k \,\ell^{T},\,$ where
the vectors  $\,k,\,\ell \in \mathbb{C}^{n+1}.\,$ Therefore,
$$
f(s,t) = x^{T}My=x^{T}k\, \ell^{T}y=k(s) \ell(t),
$$
where $\,k(s)\,$ and $\,\ell(t)\,$ are polynomials. Since,
$\,f(0,0)=1=k(0) \ell(0),\,$ and
\begin{eqnarray*}
& & det(I-sA) = f(s,0)=k(s) \ell(0) ,\\
& & det(I-tB) = f(0,t)=k(0) \ell(t)
\end{eqnarray*}
clearly,
$$
f(s,t)=k(s)\ell(0) k(0) \ell(t)=det(I-sA) \,det(I-tB).
$$
\qed

\begin{exam}
{\rm Let the matrices
$$
A=\left[
\begin{array}{ccc}
  0 & 0 & 0 \\
  0 & 1-\gamma & 1 \\
  0 & 0 & 1-\gamma \\
\end{array}%
\right] ,\;\;\;\;\;B=\left[
\begin{array}{ccc}
  0 & \gamma & 0 \\
  1/\gamma & 0 & 0 \\
  0 & 0 & 0 \\
\end{array}%
\right].$$ }
\end{exam}
We have
\begin{eqnarray*}
f(s,t)=det(I-sA-tB)& = &1+2(\gamma-1)s+(\gamma
-1)^{2}s^{2}-t^{2}+(1-\gamma)t^{2}s \\
& = & x^{T} \left[%
\begin{array}{cccc}
  1 & 0 & 1 & 0 \\
  2(\gamma -1) & 0 & 1-\gamma & 0 \\
  \gamma-1 & 0 & 0 & 0 \\
  0 & 0 & 0 & 0 \\
\end{array}%
\right] y
\end{eqnarray*}
and
$$det(I-sA)=\left(1+(\gamma-1)s\right)^{2},\;\;
\,det(I-tB)=1-t^{2}.$$ By the criterion (Proposition \ref{p3})
easily we recognize that $\,A,\,B\,$ have the CS property only for
$\,\gamma =1.\,$
\begin{rema}
~ \rm{  In equation (\ref{eqc8}), if $\,b^{T}\,a = 0\,$ then
$\,M^{2}=0, \,$ and $\,M \left( \displaystyle{\frac{ \,1\,}{\,\| b
\, \|^{2}}}\,b \right) = a.\,$ Therefore,
$$
M=P\left[%
\begin{array}{ccc}
  0 & \cdots & 1 \\
  \vdots &  & \vdots \\
  0 & \cdots & 0 \\
\end{array}%
\right]P^{-1}= P  \left[%
\begin{array}{c}
  1 \\
  0 \\
  \vdots \\
  0 \\
\end{array}%
\right]
\left[%
\begin{array}{ccc}
  0 & \cdots & 1 \\
\end{array}%
\right] P^{-1}
$$
where $\,P=
\left[%
\begin{array}{ccccc}
  a & p_2 & \cdots  & p_{n-1} & \displaystyle{\frac{ \,1\,}{\| b \,
\|^{2}}}b   \\
\end{array}%
\right]\,$ and $\,p_k, \ldots, p_{n-1}\,$ is an orthonormal basis
of $\,span\{ a, b\}^{\perp}.\,$ Then $\,P^{-1}=\left[%
\begin{array}{ccccc}
  \displaystyle{\frac{ \,1\,}{\| a \,
\|^{2}}}a & p_2 & \cdots  & p_{n-1} & b   \\
\end{array}%
\right]^{T}.\,$ }
\end{rema}

\vspace{0,4cm}

Following we note by
$\,M\left(%
\begin{array}{c}
  a_{i_1, \ldots, \,i_{p}} \\
   b_{j_1, \ldots,\, j_{q}} \\
\end{array}%
\right)\, $ the leading principal minor of order $\,p+q \,( \leq
n),\,$ which is defined by the $\,i_1, \ldots, i_{p} \,$ rows of
$\,A\,$ and $\,j_1, \ldots, j_{q}\,$ rows of $\,B,\,$ i.e.,
\begin{eqnarray*}
M\left(%
\begin{array}{c}
  a_{i_1, \ldots, \,i_{p}} \\
   b_{j_1, \ldots,\, j_{q}} \\
\end{array}%
\right)= \left|%
\begin{array}{cccccccc}
  a_{i_1 i_1} & a_{i_1 i_2}  & a_{i_1 j_1} & a_{i_1 i_3}  & \cdots   & a_{i_1 j_q}  & \cdots  & a_{i_1 i_p}  \\
  a_{i_2 i_1}  & a_{i_2 i_2}  & a_{i_2 j_1}  & a_{i_2 i_3}  & \cdots & a_{i_2 j_q}  & \cdots & a_{i_2 i_p}  \\
  b_{j_1 i_1}  &   b_{j_1 i_2}  &   b_{j_1 j_1}  &   b_{j_1 i_3}  & \cdots &   b_{j_1 j_q}  & \cdots &   b_{j_1 i_p}  \\
  a_{i_3 i_1}  & a_{i_3 i_2}  & a_{i_3 j_1}  & a_{i_3 i_3}  &  &  &  & \vdots \\
  \vdots & \vdots  &  &  & \ddots &  &  &  \\
    b_{j_q i_1}  &   b_{j_q i_2}  &  &  &  &   b_{j_q j_q}  &  &  \\
  \vdots & \vdots &  &  &  &  & \ddots &  \\
  a_{i_p i_1}  & a_{i_p i_2}  &  & \cdots &  &  &  & a_{i_p i_p}  \\
\end{array}%
\right|
\end{eqnarray*}
for $\,i_1< i_2 < j_1 < i_3 <  \cdots < j_q < \cdots < i_{p}.\,$
Thus, we clarify a determinental expression of coefficients
$\,m_{pq}\,$ in (\ref{eqc7}):
\begin{eqnarray}
\label{eqc10} m_{pq} = (-1)^{p+q} \sum_{1\leq \,i_1< \,j_1 <
\cdots
< \,j_q <\, i_{p}\leq \,n} M\left(%
\begin{array}{c}
  a_{i_1, \ldots, \,i_{p}} \\
   b_{j_1, \ldots,\, j_{q}} \\
\end{array}%
\right), \;\;\;m_{00}=1.
\end{eqnarray}

\vspace{0.3cm}

 For example, for $\, n \times n\,$ matrices $\,A\,$
and $\,B\,$ the coefficients of $\,t,\,st,\,s^{2}\,$ and
$\,s^{2}t\,$ are respectively equal to
\begin{eqnarray*}
&& m_{01}= - \sum_{1 \leq j \leq n} M(b_{j})=- \left( b_{11} +
b_{22} + \cdots + b_{nn} \,\right)= -tr B \\
&& m_{11}=
 \sum_{1\leq i < j \leq \,n} M\left(%
\begin{array}{c}
  a_{i} \\
   b_{j} \\
\end{array}%
\right) = \sum_{ {\small
\begin{array}{c} i,j=1 \\ i < j \end{array} }}^{n} \,\left(  \left|
\begin{array}{cc}
a_{ii} & a_{ij} \\
b_{ji} & b_{jj} \end{array} \right| + \left|
\begin{array}{cc}
b_{ii} & b_{ij} \\
a_{ji} & a_{jj} \end{array} \right| \right) \\
&& m_{20}=\sum_{1 \leq i,j \leq n} M(a_{ij}) = \sum_{ {\small
\begin{array}{c} i,j=1 \\ i < j \end{array} }}^{n} \, \left|
\begin{array}{cc}
a_{ii} & a_{ij} \\
a_{ji} & a_{jj} \end{array} \right|
\end{eqnarray*}
and
\begin{eqnarray*}
m_{21}=-
\sum_{1\leq i \leq \, j  \leq \,k \leq \,n} M\left(%
\begin{array}{c}
  a_{i,j} \\
   b_{k} \\
\end{array}%
\right) =- \sum_{i\leq \, j \leq \, k \leq \, n} \left(  \left|
\begin{array}{ccc}
a_{ii} & a_{ij}& a_{ik} \\
a_{ji} & a_{jj}& a_{jk} \\
b_{ki} & b_{kj} & b_{kk} \end{array} \right| + \left|
\begin{array}{ccc}
a_{ii} & a_{ij}& a_{ik} \\
b_{ji} & b_{jj} & b_{jk} \\
a_{ki} & a_{kj} & a_{kk} \end{array} \right|  + \left|
\begin{array}{ccc}
b_{ii} & b_{ij}& b_{ik} \\
a_{ji} & a_{jj} & a_{jk} \\
a_{ki} & a_{kj} & a_{kk} \end{array} \right|\right).
\end{eqnarray*}
Hence, for the matrix $\,M\,$ in (\ref{eqc7}) we have:
\newline
$M=$
\begin{eqnarray*}
{\small
~~\left[%
\begin{array}{cccccc}
  1 & -\sum M(b_{j})\;\;\;\;\;\;\;\;\;\;\;\; & \sum M(b_{i,j})\;\;\;\;\;\;\;\;\;  & \cdots & (-1)^{n-1}\sum M(b_{j_1, \ldots , j_{n-1}})\;\;  &   (-1)^{n}|B| \\
  \\
  -\sum M(a_i)\;\;\;\;\;\;\; & \sum M\left(%
\begin{array}{c}
  a_{i} \\
   b_{j} \\
\end{array}%
\right) & -\sum M\left(%
\begin{array}{c}
  a_{i} \\
   b_{j_1, j_2} \\
\end{array}%
\right) & \cdots & (-1)^{n}\sum M\left(%
\begin{array}{c}
  a_{i} \\
   b_{j_1, \ldots , j_{n-1}} \\
\end{array}%
\right) & 0 \\
\\
  \sum M \left( a_{i_1, i_2} \right) & -\sum M\left(%
\begin{array}{c}
  a_{i_1,i_2} \\
   b_{j_1} \\
\end{array}%
\right) & \vdots &  & \vdots & \vdots \\
\\
\vdots & \vdots& & & & \\
\\
  \vdots & (-1)^{n}\sum M\left(%
\begin{array}{c}
  a_{i_1,\ldots , i_{n-1}} \\
   b_{j} \\
\end{array}%
\right) & 0 &  \cdots    &  & 0 \\
\\
  (-1)^{n}|A| & 0 & 0 & \cdots &   & 0 \\
\end{array}%
\right] ~~}
\end{eqnarray*}

\vspace{0.3cm} The zeros in $\,M\,$ correspond to the coefficients
of monomials of $\,f(s,t)\,$ with degree $\,\geq n+1.\,$ These
terms are not presented in $\,det(I-sA-tB),\,$ since by
(\ref{eqc10}) the order of principal minors is greater than
$\,n.\,$ Moreover, the dimension of $\,M\,$ in (\ref{eqc7}) should
be less than $\,n+1,\,$ since the CS equation make sense for
singular matrices. \newline

Using the criterion in Proposition \ref{p3} in the above
formulation of $\,M,\,$ it is clear the next necessary and
sufficient conditions.

\vspace{0.3cm}
\begin{prop}
\label{p4}~ The $\,n \times n\,$ matrices $\,A\,$ and $\,B \,$
have the CS property if and only if
\begin{eqnarray}
\label{eqc11} \sum M (a_{i_1, \ldots , i_p}) \, \sum M (b_{j_1,
\ldots , j_q}) = \sum M\left(%
\begin{array}{c}
  a_{i_1,\ldots , i_{p}} \\
   b_{j_1, \ldots , j_{q}} \\
\end{array}%
\right),    \,\,\;\;\;\;\mbox{ for } \;\,\;p+q \leq n,~~~~~~~~~~~~
\nonumber\\
\mbox{ and ~~~~~~~~~~~~~~~~~~~~~~~~}  ~~~~~~~~~~~~~~~~~~~~~~~~~~~~~~~~~~~~~  ~~~~~~~~~~~~~~~~~~~~~~~~~~~~~~~~~~~~~~~~~~~~~~ \\
\sum M (a_{i_1, \ldots , i_p}) \, \sum M (b_{j_1, \ldots , j_q}) =
0 , \,\,\;\;\;\;\mbox{ for } \;\,\;p+q > n.
~~~~~~~~~~~~~~~~~~~~~~\nonumber
\end{eqnarray}
\end{prop}

\vspace{0.3cm}

\begin{exam}
{\rm In (\ref{eq1}) let  $\,A\,$ be a nilpotent matrix. Then,
$$
\sum M(a_i) = \sum M(a_{i,j}) = \cdots = |A| =0,
$$
 }
\end{exam}
and by Proposition \ref{p4} clearly
$$
\sum M\left(%
\begin{array}{c}
  a_{i_1,\ldots , i_{p}} \\
   b_{j_1, \ldots , j_{q}} \\
\end{array}%
\right) =0 \,\,\;\;\;\;;\;\;\;p,q=1,2, \ldots, n-1.
$$
In this case, $\,M =
\left[%
\begin{array}{c}
  1 \\
  0 \\
  0 \\
\end{array}%
\right] \,
\left[%
\begin{array}{ccccc}
  1 & b_{n-1} & \cdots  & b_1 & b_0 \\
\end{array}%
\right]. \,$

\vspace{0.4cm}

The equations (\ref{eqc11}) give also an answer to the problem
"{\it For the $\,n \times n\,$ matrix $\,A,\,$ clarify the set}
$\, CS(A)= \{ \,B \,: \;\; \;\;A\; \mbox{ and } \;B\; \mbox{
follow the CS property}\,\}. \,$

\vspace{0.4cm}

If $\, a(s) =det(I-sA) \,$ and  $\, b(t) =det(I-tB), \,$ easily we
turn out the $\,\mu$-th order derivative of polynomials at the
origin
$$
\frac{\,1\,}{p!}\,a^{(p)}(0)= \sum M (a_{i_1,  \ldots , i_p})\,,
\;\;\;\;\;\;\; \frac{\,1\,}{q!}\,b^{(t)}(0)= \sum M (b_{j_1,
\ldots , j_q}),\;\;
$$
and even
$$
\frac{\,1\,}{p!q!}\,\frac{ \partial^{\,p+q}f(0,0) }{\partial s^{p} \, \partial t^{q}} = \sum M\left(%
\begin{array}{c}
  a_{i_1,\ldots , i_{p}} \\
   b_{j_1, \ldots , j_{q}} \\
\end{array}%
\right) .
$$
Thus, if we use the Taylor's expansion of polynomials in
(\ref{eq1}), by the relationships
\begin{eqnarray*}
&& a^{(p)} (0)\, b^{(q)} (0) =  \frac{
\partial^{ p+q} f(0,0) }{\partial s^{p}\,
\partial t^{q}},\,\;\;\;\mbox{ for }\;\;\;p+q \leq n,
\\
\\
&& a^{(p)} (0)\, b^{(q)} (0) = 0,\,\;\;\;\mbox{for }\;\;\;p+ q >
n,
\end{eqnarray*}
the equations (\ref{eqc11}) arise again.

\end{document}